\documentclass[10pt,a4paper,reqno]{amsart}
\usepackage{amsfonts,amsthm,latexsym,amsmath,amssymb,amscd,epsfig}
\usepackage{graphics,graphicx}
\usepackage[all]{xy}
\ifx\undefined\xspace \usepackage{xspace} \fi

\newcommand{\nc}{\newcommand}
\renewcommand{\frak}{\mathfrak}
\providecommand{\cal}{\mathcal}
\renewcommand{\bold}{\mathbf}
\numberwithin{equation}{section}

\newcommand{\pfname}{Proof.}

                      {$\square$\vskip\medskipamount\par}

                      {$\square$\vskip\medskipamount\par}

\newtheorem{thm}{Theorem}[section]

\newtheorem{thmnn}{Theorem}

\newtheorem{corollary}[thm]{Corollary}

\newtheorem{proposition}[thm]{Proposition}

\newtheorem{lemma}[thm]{Lemma} 

\theoremstyle{definition}
\newtheorem{definition}[thm]{Definition}

\theoremstyle{definition}

\newtheorem{example}[thm]{Example} 

\nc{\Theorem}[1]{Theorem~{#1}}
\nc{\Th}[1]{({\sl Th.}~#1)}
\nc{\Thd}[2]{({\sl Th.}~{#1} {#2})}
\nc{\Theorems}[2]{Theorems~{#1} and ~{#2}}
\nc{\Thms}[2]{({\it Thms. ~{#1} and ~{#2}})}
\nc{\Lemmas}[2]{Lemma~{#1} and ~{#2}}

\nc{\manga}[6]{({\it Thms. ~ #1, ~ #2, ~ #3,\\ ~ #4, ~ #5, ~ #6})}

\nc{\Prop}[1]{({\sl Prop.}~{#1})}
\nc{\Proposition}[1]{Proposition~{#1}}
\nc{\Propositions}[2]{Propositions~{#1} and ~{#2}}
\nc{\Props}[2]{({\sl Props.}~{#1} and ~{#1})}
\nc{\Cor}[1]{({\sl Cor.}~{#1})}
\nc{\Corollary}[1]{Corollary~{#1}}
\nc{\Corollaries}[2]{Corollaries~{#1} and ~{#2}}
\nc{\Definition}[1]{Definition~{#1}}
\nc{\Defn}[1]{({\sl Def.}~{#1})}
\nc{\Lemma}[1]{Lemma~{#1}} 
\nc{\Lem}[1]{({\sl Lem.} ~{#1})} 
\nc{\Eq}[1]{equation~({#1})}
\nc{\Equation}[1]{Equation~({#1})}
\nc{\Section}[1]{Section~{#1}}
\nc{\Sections}[1]{Sections~{#1}}
\nc{\Sec}[1]{({\sl Sec.} ~{#1})} 
\nc{\Chapter}[1]{Chapter~{#1}}
\nc{\Chapt}[1]{({\sl Ch.}~{#1})}

\nc{\Ex}[1]{{\sl Ex.}~{#1}}
\nc{\Exa}[1]{{\sl Example}~{#1}}
\nc{\Example}[1]{{\sl Example}~{#1}}
\nc{\Examples}[1]{{\sl Examples}~{#1}}
\nc{\Exercise}[1]{{\sl Exercise}~{#1}}

\nc{\Rem}[1]{({\sl Rem.}~{#1})}
\nc{\Remark}[1]{{\sl Remark}~{#1}}
\nc{\Remarks}[1]{{\sl Remarks}~{#1}}
\nc{\Note}[1]{{\sl Note}~{#1}}

\nc{\Conjecture}[1]{Conjecture~{#1}}
\nc{\Claim}[1]{Claim~{#1}}

\nc \Proof{{  \it Proof. }}

\nc{\xmu}{\mu}
\nc{\w}{\omega}

\nc \Ab{{\ensuremath{\bold A}}}
\nc \ab{{\ensuremath{\bold a}}}
\nc \bb{{\ensuremath{\bold b}}}
\nc \cb{{\ensuremath{\bold c}}}
\nc \Bb{{\ensuremath{\bold B}}}
\nc \Gb{{\ensuremath{\bold G}}}
\nc \Qb{{\ensuremath{\bold Q}}}
\nc \Rb{{\ensuremath{\bold R}}} \nc \Cb{{\ensuremath{\bold C}}} 
\nc \Eb{{\ensuremath{\bold E}}}
\nc \eb{{\ensuremath{\bold e}}}
\nc \Db{{\ensuremath{\bold D}}}
\nc \Fb{{\ensuremath{\bold F}}}
\nc \ib{{\ensuremath{\bold i}}}
\nc \jb{{\ensuremath{\bold j}}}
\nc \kb{{\ensuremath{\bold k}}}
\nc \nb{{\ensuremath{\bold n}}}
\nc \rb{{\ensuremath{\bold r}}}
\nc \Pb{{\ensuremath{\bold P}}}
\nc \pb{{\ensuremath{\bold p}}}
\nc \SPb{{\ensuremath{\bold {SP}}}}
\nc \Zb{{\ensuremath{\bold Z}}} 
\nc \zb{{\ensuremath{\bold z}}} 
\nc \gb{{\ensuremath{\bold g}}} 
\nc \fb{{\ensuremath{\bold f}}} 
\nc \ub{{\ensuremath{\bold u}}} 
\nc \vb{{\ensuremath{\bold v}}} 
\nc \yb{{\ensuremath{\bold y}}} 
\nc \xb{{\ensuremath{\bold x}}} 
\nc \xib{{\ensuremath{\bold \xi}}} 
\nc \Nb{{\ensuremath{\bold N}}} 
\nc \Hb{{\ensuremath{\bold H}}} 
\nc \wb{{\ensuremath{\bold w}}} 
\nc \Wb{{\ensuremath{\bold W}}} 
\nc \syz{{\mathbf {syz}}}
\nc \bnoll{{\ensuremath{\bold 0}}} 

\nc \mf{\frak m} \nc \mh{\hat{\m}} 
\nc \nf{\frak n}
\nc \Of{\frak O}
\nc \rf{\frak r}
\nc \mufr{{\mathbf \mu}}
\nc \hf{\frak h} 
\nc \qf{\frak q} 
\nc \bfr{\frak b} 
\nc \kfr{\frak k} 
\nc \pfr{\frak p} 
\nc \af{\frak a }
\nc \cf{\frak c }
\nc \sfr{\frak s} 
\nc \ufr{\frak u} 
\nc \g{\frak g} 
\nc \gA{\g_{\Ao}} 
\nc \lfr{\frak l}
\nc \afr{\frak a}
\nc \gfh{\hat {\frak g}}
\nc \gl{\frak { gl }}
\nc \Sl{\frak {sl}}
\nc \SU{\frak {SU}}
\nc{\Homf}{\frak{Hom}}

\newcommand{\on}{\operatorname}
\nc\hankel{\on {Hankel}}
\nc\row{\on {row\ }}
\nc\nullity{\on {nullity }}
\nc\col{\on {col\ }}
\nc\rowm{\on {Row \ }}
\nc\loc{\on {lc \ }}
\nc\nullo{\on {null\ }}
\nc\Nul{\on {Nul\ }}
\nc \Ann {\on {Ann }}
\nc \Ass {\on {Ass \ }}
\nc \Coker {\on {Coker}}
\nc \Co{\on C}
\nc \Homo{\on {Hom}}
\nc \Ker {\on {Ker}}
\nc \omod{\on {mod}}
\nc \No {\on N}
\nc \NN {\on {NN}}
\nc \NGo {\on {NG}}
\nc \Oo {\on O}
\nc \ch {\on {ch}}
\nc \rko {\on {rk}}
\nc \Sing {\on {Sing\ }}
\nc \Reg {\on {Reg}}
\nc \CoI {\on {CI}}
\nc \CoM {\on {CM}}
\nc \Gor {\on {Gor}}
\nc \Type {\on {Type}}
\nc \can {\on {can}}
\nc \Top {\on {T}}
\nc \Tr {\on {Tr}}
\nc \rel {\on {rel}}
\nc \tr {\on {tr}}
\nc \sgn {\on {sgn }}
\nc \trdeg {\on {tr.deg}}
\nc \codim {\on {codim }}
\nc \coht {\on {coht}}
\nc \divo {\on {div \ }}
\nc \coh {\on {coh}}
\nc \Clo {\on {Cl}}
\nc \embdim{\on {embdim}}
\nc \embcodim{\on {embcodim \ }}
\nc \qcoh {\on {qcoh}}
\nc \grad {\on {grad}\ }
\nc \grade {\on {grade}}
\nc \hto {\on {ht}}
\nc \depth {\on {depth}}
\nc \prof {\on {prof}}
\nc \reso{\on {res}}
\nc \ind{\on {ind}}
\nc \prodo{\on {prod}}
\nc \coind{\on {coind}}
\nc \Con{\on {Con}}
\nc \Crit{\on {Crit}}
\nc \Der{\on {Der}}
\nc \Char{\on {Char}}
\nc \Ch{\on {Ch}}

\nc \Ext{\on {Ext}}
\nc \Eo{\on {E}}
\nc \End{\on {End}}
\nc \ad{\on {ad}}
\nc \Ad{\on {Ad}}
\nc \gr{\on {gr}}
\nc \Fo{\on {F}}
\nc \Gr{\on {Gr}}
\nc \Go{\on {G}}
\nc \GFo{\on {GF}}
\nc \Glo{\on {Gl}}
\nc \Ho{\on {H}}
\nc \CMo{\on {\CM}}
\nc \SCM{\on {SCM}}
\nc \hol{\on {hol}}
\nc{\sgd}{\on{sgd}}
\nc \supp{\on {supp}}
\nc \ssupp{\on {s-supp}}
\nc \singsupp{\on {singsupp}}
\nc \msupp{\on {msupp}}
\nc \spec{\on {spec}}
\nc \spano{\on {span }}
\nc \Span{\on {Span }}
\nc \Max{\on {Max}}
\nc \Min{\on {Min}}
\nc \Mod{\on {Mod}}
\nc \Rad {\on {Rad}}
\nc \rad {\on {rad}}
\nc \rank {\on {rank}}
\nc \range {\on {range}}
\nc \Slo{\on {SL}}
\nc \soc {\on {soc}}
\nc \Irr {\on {Irr}}
\nc \Imo {\on {Im}}
\nc \SSo{\on {SS}}
\nc \lub{\on {lub}}
\nc \gldim{\on {gl.d.}}
\nc \pdo{\on {p.d.}} 
\nc \ido{\on {i.d.}} 
\nc \dSSo{\dot {\SSo}}
\nc \So{\on S}
\nc \Io{\on I}
\nc \Jo{\on J}
\nc \jo{\on j}
\nc \Ko{\on K}
\nc \PBW{\Ac_{PBW}}
\nc \Ro{\on R}
\nc \To{\on T}
\nc \Ao{\on A}

\nc \Do{{\on D}}
\nc \Bo{\on B}
\nc \Po{\on P}
\nc \Qo{\on Q}
\nc \Zo{\on Z}
\nc \U{\on U}
\nc \wt{\on {wt}}
\nc \Uh{\hat {\U}}
\nc \T{\on T}
\nc \Lo{\on L}
\nc{\dop}{\on d}
\nc{\eo}{\on e}
\nc{\ado}{\on{ad}}
\nc{\Tot}{\on{Tot}}
\nc{\Aut}{\on{Aut}}
\nc{\sinc}{\on {sinc}}

\nc{\overrightleftarrows}[2]{\overset{#1}{\underset{#2}{\rightleftarrows}}}

\nc{\CCF}{\cal{CF}}
\nc{\CDF}{\cal{DF}}
\nc{\CHC}{\check{\cal C}}

\nc{\Cone}{\on{Cone}}
\nc{\dec}{\on{dec}}
\nc{\Diff}{\on{Diff}}
\nc{\dirlim}{\underset{\to}{\on{lim}}}
\nc{\dpar}{\partial}
\nc{\GL}{\on{GL}}
\nc{\CGr}{\cal{G}r}
\nc{\pr}{\on{pr}}
\nc{\semid}{|\!\!\!\times}
\nc{\Hom}{\on{Hom}}
\nc \RHom{\on {RHom}}

\nc \Proj{\mathrm {Proj\ }}
\nc \proj{\mathrm {proj}}
\nc{\Id}{\on{Id}}
\nc{\id}{\on{id}}
\nc{\Ima}{\on{Im}}
\nc{\invtimes}{\underset{\gets}{\otimes}}
\nc{\invlim}{\underset{\gets}{\on{lim}}}
\nc{\Lie}{\on{Lie}}
\nc{\re}{\on{Re }}
\nc{\Pic}{\on{Pic }}
\nc{\LPic}{\on{LPic }}
\nc{\Sch}{\on{Sch}}
\nc{\Sh}{\on{Sh}}
\nc{\Set}{\on{Set}}
\nc{\spo}{\on{sp\  }}
\nc{\Spec}{\on{Spec}}
\nc{\mSpec}{\on{mSpec}}
\nc{\Specb}{\bold {Spec}}
\nc{\Projb}{\bold {Proj}}
\nc{\Specan}{\on{Specan}}
\nc{\Spo}{\on{Sp}}
\nc{\Spf}{\on{Spf}}
\nc{\sym}{\on{sym}}
\nc{\symm}{\on{symm}}
\nc{\rop}{\on{r}}
\nc{\Td}{\on{Td}}
\nc{\Tor}{\on{Tor}}

\nc{\Artin}{\cal{A}rtin}
\nc{\Dgcoalg}{\cal{D}gcoalg}
\nc{\Dglie}{\cal{D}glie}
\nc{\Ens}{\cal{E}ns}
\nc{\Fsch}{\cal{F}sch}
\nc{\Groupoids}{\cal{G}roupoids}
\nc{\Holie}{\cal{H}olie}
\nc{\Mor}{\cal{M}or}

\nc{\CF}{\ensuremath{\cal{F}}}
\nc \Kc{\ensuremath{\cal K}}
\nc \Lc{{\ensuremath{\cal L}}}
\nc \lcc{{\mathcal l}} 
\nc \CC{{\ensuremath{\cal C}}} 
\nc \Cc{{\ensuremath {\cal C}}}
\nc \Pc{{\ensuremath{\cal P}}}
\nc \Dc{\ensuremath{\mathcal D}}
\nc \Ac{{\ensuremath{\cal A}}} 
\nc \Bc{{\ensuremath{\cal B}}}
\nc \Ec{{\ensuremath{\cal E}}}
\nc \Fc{{\ensuremath{\cal F}}}
\nc \Mcc{{\ensuremath{\cal M}}} 
\nc \hM{\hat{\Mcc}} 
\nc \bM{\bar {\Mcc}} 
\nc\hbM{\hat{\bar \Mcc}}  
\nc \Nc{{\ensuremath{\cal N}}}
\nc \Hc{{\ensuremath{\cal H}}} 
\nc \Ic{{\ensuremath{\cal I}}} 
\nc \Oc{\ensuremath{{\cal O}}}
\nc \Och{\hat{\cal O}} 
\nc \Sc{{\ensuremath{{\cal S}}}}
\nc \Tc{\ensuremath{{\cal T}}} 
\nc \Vc{{\ensuremath{{\cal V}}}} 
\nc{\CA}{{\ensuremath{{\cal A}}}}
\nc{\CB}{{\ensuremath{{\cal B}}}}
\nc{\C}{{\ensuremath{{\cal F}}}}
\nc{\Gc}{{\ensuremath{{\cal G}}}}
\nc{\CH}{\ensuremath{\mathcal H}}
\nc{\CI}{{\ensuremath{{\cal I}}}}
\nc{\CM}{{\ensuremath{{\cal M}}}}
\nc{\CN}{{\ensuremath{{\cal N}}}}
\nc{\CO}{{\ensuremath{{\cal O}}}}
\nc{\Rc}{{\ensuremath{{\cal R}}}}
\nc{\CT}{{\ensuremath{\mathcal T}}}
\nc{\CU}{\ensuremath{{\cal U}}}
\nc{\CV}{\ensuremath{{\cal V}}}
\nc{\CZ}{\ensuremath{{\cal Z}}}
\nc{\Homc}{\ensuremath{{\cal {Hom}}}}

\nc{\fa}{\frak{a}}
\nc{\fA}{\frak{A}}
\nc{\fg}{\frak{g}}
\nc{\fh}{\frak{h}}
\nc{\fI}{\frak{I}}
\nc{\fK}{\frak{K}}
\nc{\fm}{\frak{m}}
\nc{\fP}{\frak{P}}
\nc{\fS}{\frak{S}}
\nc{\ft}{\frak{t}}
\nc{\fX}{\frak{X}}
\nc{\fY}{\frak{Y}}

\nc{\bF}{\bar{F}}
\nc{\bCP}{\bar{\cal{P}}}
\nc{\bm}{\mbox{\bf{m}}}
\nc{\bT}{\mbox{\bf{T}}}
\nc{\hB}{\hat{B}}
\nc{\hC}{\hat{C}}
\nc{\hP}{\hat{P}}
\nc{\htest}{\hat P}

\nc{\nen}{\newenvironment}
\nc{\ol}{\overline}
\nc{\ul}{\underline}
\nc{\ra}{\to}
\nc{\lla}{\longleftarrow}
\nc{\lra}{\longrightarrow}
\nc{\Lra}{\Longrightarrow}
\nc{\Lla}{\Longleftarrow}
\nc{\Llra}{\Longleftrightarrow}
\nc{\hra}{\hookrightarrow}
\nc{\iso}{\overset{\sim}{\lra}}

\nc{\dsize}{\displaystyle}
\nc{\sst}{\scriptstyle}
\nc{\tsize}{\textstyle}
\nen{exa}[1]{\label{#1}{\bf Example.\ } }{}

\nen{rem}[1]{\label{#1}{\em Remark.\ } }{}
\nen{exer}[1]{\label{#1}{\em Exercise.\ } }{}

\newcommand {\bC} {\mathbb C}

\newcommand {\bZ} {\mathbb Z}

\begin{document}
\numberwithin{equation}{section}
\title[on decomposition factors of D-modules]{Decomposition factors of D-modules on hyperplane configurations in general position}

\author{Tilahun Abebaw and Rikard B{\o}gvad}
\address{Department of Mathematics, Addis Ababa University and Stockholm University}
\email{tabebaw@math.aau.edu.et, abebaw@math.su.se,}
\address{Department of Mathematics, Stockholm University}\email{rikard@math.su.se}
\date{}
\maketitle
\begin{abstract} Let $\alpha_{1},...,\alpha_{m}$ be linear functions
on $\bC^{n}$ and
$\rm{X=\bC^{n}\setminus V(\alpha)},$ where $\alpha=\prod\limits_{i=1}^{m}\alpha_{i}$ and 
$\rm{V(\alpha)=\{p\in\bC^{n}:\alpha(p)=0\}}$. The coordinate ring
$\rm{\mathcal{O}_{X}}=\bC[x]_{\alpha}$ of $\rm{X}$ is a holonomic $A_{n}$-module,
where $A_{n}$ is the n-th Weyl algebra and since holonomic
$A_{n}$-modules have finite length, $\rm{\mathcal{O}_{X}}$ has finite
length. We consider a ''twisted'' variant of this $A_{n}$-module
which is also holonomic. Define $\rm{M_{\alpha}^{\beta}}$ to be the
free rank 1 $\bC[x]_{\alpha}$-module on the generator
$\alpha^{\beta}$ (thought of as a multivalued function), where
$\alpha^{\beta}=\alpha_{1}^{\beta_{1}}...\alpha_{m}^{\beta_{m}}$ and
the multi-index $\beta=(\beta_{1},...,\beta_{m})\in\bC^{m}$. 
It is straightforward to describe the decomposition
factors of $\rm{M_{\alpha}^{\beta}}$,  when the linear functions $\alpha_{1},...,\alpha_{m}$ define a normal crossing hyperplane configuration, and we use this to give a sufficient
 criterion on $\beta$ for the irreducibility of $\rm{M_{\alpha}^{\beta}}$, in terms of numerical data for a resolution of the singularities of $V(\alpha).$
\end{abstract}

\section{Introduction}

\begin{definition}\label{module}Let $\rm{V}\cong \bC^{n}$ be a finite-dimensional complex vector space.
\begin{itemize}
\item[(i)]A (finite) hyperplane configuration  $\rm{A=\{H_{1},...,H_{m}\}}$ in $V$ is a finite set of affine hyperplanes $H_i=V(\alpha_i)$ in $V$, where $\alpha_i : V\to \bC$ is a non-zero polynomial of degree one.
\item[(ii)]A hyperplane configuration $\rm{A}$ in $V$ is said to be in {\it general position}, if
$$\{H_{i_{1}}, . . . , H_{i_{p}}\} \subset A,\  p\leq n \implies {\rm  dim}_{\bC}(H_{i_{1}}\cap...\cap H_{i_{p}}) = n - p$$ and 
$$\{H_{i_{1}}, . . . , H_{i_{p}}\} \subset A, p>n \implies H_{i_{1}}\cap...\cap H_{i_{p}} = \emptyset.$$ 
\item[(iii)] A hyperplane configuration $\rm{A}$ is said to be a {\it normal crossing configuration}, if for any set of hyperplanes $\rm{H_{i_{1}},...,H_{i_{k}}}$ that contain $p$, there is an affine change of local coordinates at $p$ to new coordinates $x_{1},...,x_{n}$, such that $\alpha_{i_{1}}=x_{1},...,\alpha_{i_{k}}=x_{k}.$ In particular any point  $p\in V$ is contained in at most $n={\rm dim} V$ of the hyperplanes in $A$.
\end{itemize}
\end{definition}

If $A$ is a normal crossing configuration, then the divisor $\alpha=\prod\limits_{i=1}^{m}\alpha_{i}$ is a normal crossings divisor(see \cite{HJ1}). Any  configuration in general position is a normal crossing configuration. The intersections $\rm{H\neq\bC^{n}}$ of a subset of the hyperplanes in $A$ are called the {\it flats} of the configuration. 

\begin{example} Let $\rm{V=\bC^{2}}.$ Then a set of lines in $\bC^{2}$ defines 
\begin{itemize}
\item[(i)] a hyperplane configuration in general position if no two are parallel
and no three meet at a point and ;
\item[(ii)] a normal crossing configuration if no three meet at a point.
\end{itemize}
\end{example}

There is an extensive literature on hyperplane arrangements {\cite{Stanley}. The problem studied in this note is the following.
Let $A=\{ H_i=V(\alpha_i),\ i=1,...,m\} $, be a hyperplane configuration, and set $\rm{X=\bC^{n}\setminus V(\alpha)},$ where $\alpha=\prod\limits_{i=1}^{m}\alpha_{i}$. The coordinate ring $\rm{\mathcal{O}_{X}}$ of $\rm{X}$ is the localization $\bC[x]_{\alpha}:=\bC[x_1,...,x_n]_{\alpha}$ and this is a holonomic $A_{n}$-module, where $A_{n}$ is the n-th Weyl algebra. Consider the following twisted variant of this $A_{n}$-module, that corresponds to the multivalued function $\alpha^\beta=\alpha_1^{\beta_1}...\alpha_m^{\beta_m}$. 
\begin{definition}
\label{Def11}The module $\rm{M_{\alpha}^{\beta}}$ is (as a $\bC[x]_{\alpha}-$module) the free rank 1 $\bC[x]_{\alpha}$-module on the generator $\alpha^{\beta}$, where $\beta=(\beta_{1},...,\beta_{m})\in\bC^{m}$. It is furthermore an $A_{n}$-module, defining 
 $$\partial_{j}(\alpha^{\beta})=\sum_{i=1}^{m} \beta_{i}\frac{\partial_{j}(\alpha_{i})}{\alpha_{i}}\alpha^{\beta}$$ for $j=1,2,...,n$ and extending to an action of $A_{n}$ on $\rm{M_{\alpha}^{\beta}}$.
\end{definition}
 The $A_{n}$-module $\rm{M_{\alpha}^{\beta}}$ is holonomic and so has finite length. The decomposition factors have support on the flats associated to the configuration (Proposition~\ref{Prop32}). In continuation of \cite{TARB}, where we treated the case $n=2$, and inspired by \cite{CDCP1} where incidentally the case of $\beta\in\bZ^m$ is treated, our main interest lies in finding the decomposition factors of $\rm{M_{\alpha}^{\beta}}$. If  $A$ is a normal crossing configuration, this is not difficult. In Theorem~\ref{Thm1},  
 we prove that in that case $\rm{M_{\alpha}^{\beta}}$ has exactly one decomposition factor with support on a flat $H=\cap_{k=1}^{r}H_{i_k}$ iff the associated $\beta_{i_k}, \ k=1,...,r$ are integers.  Normal crossing configurations occur when the singularities of $V(\alpha)$ are resolved. Applied to a such a resolution,  Theorem~\ref{Thm1} together with the decomposition theorem for D-modules,  then gives our main result Theorem \ref{Thm:irr}: a sufficient criterion for $\rm{M_{\alpha}^{\beta}}$ to be irreducible for an arbitrary hyperplane arrangement, in terms of numerical data of the resolution and $\beta$.

\section{Support of D-modules}\label{s-2}

\subsection{Definition}Let X be a smooth algebraic variety. We  denote by $\rm{\mathcal{D}_{X}}$ the sheaf of differential operators on X.  If  ${\rm{X}}=\bC^{n}$ this is the same as the sheaf on $\bC^{n}$ associated to the Weyl algebra $A_{n}$. If $X$ is an affine variety, we will by abuse of notation, identify $\mathcal{O}_{X}$ with $\Gamma(X, \mathcal{O}_{X})$ and an $\mathcal{O}_{X}$-module sheaf with its global sections as a module over $\Gamma(X, \mathcal{O}_{X})$. If X is an affine open subset of $\bC^{n}$ defined by $0\neq \alpha\in\bC[x],$ then  $\rm{\mathcal{O}_{X}=\bC[x]_{\alpha}}$ and $\rm{\mathcal{D}_{X}=\bC[x]_{\alpha}\otimes_{\bC[x]}A_{n}}$.

 A holonomic $\rm{\mathcal{D}_{X}}$-module $M$ is in particular an $\rm{\mathcal{O}_{X}}$-module and as such may be shown to have support on a closed variety $Z={\rm Supp}M$(see \cite{BJE1}). If $\mathcal{I}\subset \mathcal{O}_X$ is the ideal associated to $Z$, any local section of $M$ is annihilated by a large enough power of $\mathcal{I}$, and $Z$ is the minimal closed subset with this property.

 We will see later 
 that the support of the irreducible factors of $\rm{M_{\alpha}^{\beta}}$ are intersections of hyperplanes. 
\begin{example}
\begin{itemize}
\item[(i)] Let $\rm{M_{1}}=\bC[x,y]_{xy}/(\bC[x,y]_{x}+\bC[x,y]_{y})$. Then  
$$\rm{SuppM_{1}=V(x,y)=(0,0)}.$$
\item[(ii)]Let $\rm{M_{2}=\bC[x,y]_{x}/\bC[x,y]}$. Then  $\rm{SuppM_{2}=V(x)=\{(0,y):y\in\bC\}}$. 
\item[(iii)]Let $\rm{M_{3}=\bC[x,y]}$. Then $\rm{SuppM_{3}=V(0)=\bC^{2}}$.
\end{itemize}
\end{example}
\subsection{Some Basic Properties}
Let $M$ be a $\rm{\mathcal{D}_{X}}$-module and $\rm{U\subset X}$ be an open subset. Then define  $\rm{M_{\vert_{U}}=:\mathcal{O}_{U}\otimes_{\mathcal{O}_{X}}M}$;  this is a $\mathcal{D}_{\rm{U}}$-module.
\begin{lemma}\label{Lem32}Let $\rm{U\subset X}$ be an open subset. 
\begin{itemize}
\item[(i)] $\rm{SuppM_{\vert_{U}}=U\cap {\rm {Supp}}M}$.
\item[(ii)]$\rm{M_{\vert_{U}}=0\Leftrightarrow SuppM\subset X-U=:Z}$.
\end{itemize}
\end{lemma}
\begin{proof} (i) is clear by the above description. Let $\mathcal{I }$ be the ideal of a closed subvariety $Z$, and let $j\colon U\to X$ be the inclusion. For any $\mathcal{O}_{X}$-module  $\mathcal{M}$ there exists an exact sequence of $\mathcal{O}_{X}$-modules
\begin{equation}
\label{eq:loccoh}
\rm{\Gamma_{\rm{Z}}M\subset M\longrightarrow j_*M_{\vert_{U}}},
\end{equation} 
where (locally) $\rm{\Gamma_{Z}M=\{m\in M:\exists r,\mathcal{I }^{r}m=0\}}$. If $\rm{M_{\vert_{U}}=0}$, then  $\rm{\Gamma_{Z}M=M}$ and this proves (ii) in one direction. The other direction is immediate. 
\end{proof} 
It is easy to see that localization preserves irreducibility of $\mathcal{D}_{X}$-modules.(Cf. also Lemmas 2.1-2 of \cite{Kashiwara2}; since the formulation there is different and without proof, we give one.)
\begin{proposition}\label{Prop21}If M is an irreducible $\mathcal{D}_{X}$-module and $\rm{U\subset X}$  an open subset, then $\rm{M_{\vert_{U}}=:\mathcal{O}_{U}\otimes_{\mathcal{O}_{X}}M}$ is an irreducible $\mathcal{D}_{\rm{U}}$-module. 
\end{proposition}
\begin{proof} Let $j\colon U\to X$ be the inclusion and denote the adjunction in  (\ref{eq:loccoh}) by $\delta\colon M\to j_*\rm{M_{\vert_{U}}}$. Since $M$ is irreducible $\delta$ is either $0$ or an injection. 
In the first case $M$ is $Z$-torsion and $M\vert_U=0$, and the lemma is trivially true. 
So assume that $\delta$ is an injection, and that $K\subset \rm{M_{\vert_{U}}}$ is an irreducible sub $\rm{\mathcal{D}_{U}}$-module. Then we have an injection $\gamma \colon j_*K\to  j_*\rm{M_{\vert_{U}}}$. 
The intersection $\delta(M)\cap \gamma(j_*K)$ is either $0$ or $\delta(M)$. In the second case $\delta(M)\subset  \gamma(j_*K)$. By restricting to $U$, and using that $\delta\vert_U$ is the identity and $\gamma\vert_U$ is the inclusion $K\subset \rm{M_{\vert_{U}}}$, we get that $M\vert_U=K$. Hence the lemma is true in this case too. The first case implies that $\gamma: j_*K\to  j_*\rm{M_{\vert_{U}}}/\delta(M)$ is an inclusion, which implies that $ j_*K$ is a $Z$-torsion module. Hence  $K=(j_*K)\vert _U=0$.
\end{proof}
Since localization is an exact functor we have the following corollary. 
We will use the notations $\rm{DF(M)}$ for the set of decomposition factors of $\rm{M}$ and $\rm{c(M)}$ for the number of  decomposition factors.
\begin{corollary}
\label{cor:open}
\begin{itemize} 
\item[(i)] $\rm{DF(M_{\vert_{U}})=\{M_{i}\in DF(M):SuppM_{i}\cap U\neq \emptyset\}}$. 
\item[(ii)] $\rm{c(M_{\vert_{U}})\leq c(M)}$.
\end{itemize}
\end{corollary}

\subsection{The support of decomposition factors of $\rm{M_{\alpha}^{\beta}}$}
That the support of decomposition factors of the modules  we study are intersections of hyperplanes is easily seen by general arguments, but it is also possible to give a direct argument, as we will now do. 
First we will do a reduction. 

Let $V_1=\cap_{i=1}^m H_i$, and assume that this affine space has positive dimension. Choose a vector space complement $V_2$ to the affine subspace $V_1$, such that $\bC^n\cong V_1\oplus V_2 $, and denote the restriction of a linear function $\alpha$ to $V_2$ by $\tilde \alpha$. 

\begin{lemma}[Reduction lemma]
\label{reductionlemma}
\begin{itemize}
\item[i)]  $\rm{M_{\alpha}^{\beta}}$ is isomorphic to the exterior tensor product $\rm{M_{\tilde \alpha}^{\beta}}\widehat{\otimes} {\mathcal {O}}_{V_1}$, as $A_n$-modules.
\item[ii)] If $N_i,\ i=1,...,r$ are the decomposition factors of $\rm{M_{\tilde \alpha}^{\beta}}$, then 
 $\rm{N_i}\widehat{\otimes} {\mathcal {O}}_{V_1},\ i=1,...,r$ are the decomposition factors of $\rm{M_{ \alpha}^{\beta}}$.
\item[iii)] $\rm{supp(N_i}\widehat{\otimes} {\mathcal {O}}_{V_1})=(\rm{supp N_i})\times V_1, \  i=1,...,r$.
\end{itemize}
\end{lemma}
\begin{proof}i) follows by an affine change of coordinates. By Corollary 2.5 of \cite{TARB} a decomposition factor of $\rm{M_{\alpha}^{\beta}}$ is the external tensor product of a decomposition factor of $\rm{M_{\tilde\alpha}^{\beta}}$ with ${\mathcal {O}}_{V_1}$.  This gives ii). Finally iii) is immediate.
\end{proof}
\begin{proposition}\label{Prop32} Decomposition factors of $\rm{M_{\alpha}^{\beta}}$ have support on intersections of hyperplanes.
\end{proposition}
\begin{proof}  Use the notation $\rm{H_{S}=\{p\in \bC^{n}:\alpha_{i}(p)=0,i\in S\}}$, where $S\subset\{1,2,...,m\}$. By the lemma we can reduce to the case when $V_1$ is a point or $V_1=\emptyset$. 
Make induction on the number $m$ of hyperplanes;  for $m=0$ the statement is immediate from the irreducibility of $\bC[x]$ as an $A_n$-module.  Assume it is true for an arbitrary arrangement with $m$ hyperplanes.  Then we have to prove the following: if N is a decomposition factor of $\rm{M_{\alpha}^{\beta}},$ where $\alpha=\prod\limits_{i=1}^{m+1}\alpha_{i}$, then $\rm{H:=SuppN}$ is an intersection of hyperplanes. 

By our reduction, there exists a hyperplane $\rm{H_{j}}$ such $\rm{H\nsubseteq H_{j}}$, except in the case where $H$ equals $\cap_{i=1}^{m+1}H_i$, and so in particular is a hyperplane intersection, and there is nothing to prove. Let $\rm{U=\bC^{n}\setminus H_{j}}$. Then by Proposition~\ref{Prop21}, $\rm{N|_{U}}$ is non-zero and so a decomposition factor of $\rm{M_{\alpha}^{\beta}|_{U}}.$ On $U$, $\alpha_{j}$ is invertible and by  Lemma~\ref{Prop33} below, $\gamma\colon \rm{M_{\alpha}^{\beta}|_{U}\cong M_{\bar{\alpha}}^{\bar{\beta}}|_{U}},$ where $ \bar{\alpha}=\prod_{i=1,i\neq j}^{m+1}\alpha_{i}$ and $\bar{\beta}=(\beta_{1},...\widehat{\beta_{j}},...,\beta_{m+1}).$ There is, again by Proposition~\ref{Prop21}, a decomposition factor $N_1$ of $M_{\bar{\alpha}}^{\bar{\beta}}$, such that $N_1\vert_U=\gamma(N\vert_ U)$. By induction,   $\rm{supp N_1}=H_S$  is an intersection of hyperplanes. Since twisting modules by automorphisms as in Lemma~\ref{Prop33}  preserves support, $\rm{supp N}\cap U=H_S\cap U$. Since $N$ is irreducible,  $\rm{supp N}$ is an irreducible subvariety, as is $H_S$, and so $\rm{supp N}=\overline{\rm{supp N}\cap U}=\overline{H_S\cap U}=H_S.$
This proves the induction step. 
\end{proof} 
\begin{definition} Suppose that $\theta$ is an automorphism of $\rm{\mathcal{D}_{X}}$. If M is a $\rm{\mathcal{D}_{X}}$-module, $\rm{\theta^{*}M}$ is defined to be the $\rm{\mathcal{D}_{X}}$-module which consists locally of the same sections as M, but on which (local sections of) $\rm{\mathcal{D}_{X}}$ acts by $\theta$: if $\rm{P\in \mathcal{D}_{X}}$, $\rm{m\in\theta^{*}M}$, then $Pm=\theta(P)m$.
\end{definition}
The following Lemma is clear.
\begin{lemma}\label{Lemm31}Let $\theta:\rm{\mathcal{D}_{X}\longrightarrow \mathcal{D}_{X}}$ denote an automorphism that is the identity on $\rm{\mathcal{O}_{X}}$. Suppose M has decomposition factors $\rm{M_{i}}$, $i=1,...,l$. Then $\rm{\theta^{*}M}$ has decomposition factors $\rm{\theta^{*}M_{i}}$, $i=1,...,l$. In particular $\rm{c(M)=c(\theta^{*}M)}$. The support of $\rm{\theta^{*}M_{i}}$ equals the support of $\rm{M_{i}}$.
\end{lemma}
We apply this in the following situation.
\begin{proposition} \label{Prop33}Suppose that $X=\bC^n$, so that $\mathcal O_X=\bC[x]$ and let $U=X-V(\alpha_{1},...,\alpha_{l})$. Then $\rm{c({M_{\alpha}^{\beta}}_{\vert U})=c({M_{\bar{\alpha}}^{\bar{\beta}}}_{\vert U})}$, where $\tilde{\alpha}=\alpha_{l+1}...\alpha_{m}$ and $\tilde{\beta}=(\beta_{l+1},...,\beta_{m})$. 
\end{proposition}
\begin{proof} By the preceding lemma, it suffices to construct an automorphism $\rm{\theta:\mathcal{D_{\rm{U}}}\rightarrow \mathcal{D_{\rm{U}}}}$  that is the identity on $\rm{\mathcal{O}_{U}}$ and  such that $\rm{\theta^{*}(M_{\beta}^{\beta}|_{U})\cong M_{\bar{\alpha}}^{\bar{\beta}}|_{U}}$. Further it is by induction enough to assume that $\rm{U=X-V(\alpha_{1})}$. Put $\alpha^{\beta}=\alpha_{1}^{\beta_{1}}\tilde{\alpha}^{\tilde{\beta}}$, where $\tilde{\alpha}=\alpha_{2}...\alpha_{m}$ and $\tilde{\beta}=(\beta_{2},...,\beta_{m}).$ Then $\rm{{M_{\alpha}^{\beta}}{\vert _U}=\bC[x]_{\alpha}\alpha_{1}^{\beta_{1}}{\tilde{\alpha}}^{\tilde{\beta}}}$ where $\alpha_{1}$ is invertible. Define $\rm{\theta:\mathcal{D_{\rm{U}}}\longrightarrow \mathcal{D_{\rm{U}}}}$ by: 
 $$\rm{\theta(D)=D+\frac{\beta_{1}D(\alpha_{1})}{\alpha_{1}}}$$ for all $\rm{D\in Der_{\bC}(\bC[x])\subset \mathcal{D_{\rm{U}}}}$ and  $\theta(r)=r$ for all $\rm{r\in \mathcal{O}_{U}}$, and
extend this to an endomorphism of $\rm{\mathcal{D_{\rm{U}}}}$. This is the desired automorphism, since it is easily checked that it has an inverse, and that
the map $$\rm{\rho:\theta^{*}({M_{\tilde{\alpha}}^{\tilde{\beta}}}{\vert _U})\longrightarrow{M_{\alpha}^{\beta}}{\vert _U}},$$ defined by $\rho(r\tilde{\alpha}^{\tilde{\beta}})=r\alpha^{\beta},$ is a $\rm{\mathcal{D_{\rm{U}}}}$-isomorphism. 
\end{proof}

\section{Decomposition factors on normal crossings}\label{s-3} 
We will now describe the decomposition factors of $\rm{M_{\alpha}^{\beta}}$ for a normal crossing hyperplane configuration.

\begin{lemma}\label{lemma:unique}
\label{blem} Let $\rm{A=\{H_{1},...,H_{m}\}}$ be a nonempty  normal crossing hyperplane configuration in $\bC^{n}$.  If  $\cap_{i=1}^m H_i\neq \emptyset$, then $m\leq n$ and the arrangement is central. Assume that $\cap_{i=1}^m H_i=\emptyset$. Then, for each intersection $\rm{H(\neq\bC^{n})}$ of a subset of the hyperplanes in $A$,
\begin{itemize}
\item[(i)]  there exists at least one hyperplane $\rm{H_{i}}$, such that $\rm{H\nsubseteq H_{i}}$, and
\item[(ii)]there is  a unique subset $S=\{ i_1,...,i_k\}\subset \{1,...,m\}$ such that  $H=H_S=\cap_{k=1}^{r}H_{i_k}$. Furthermore  $H_j\supset H\iff j\in S$.
\end{itemize}
\end{lemma}
The following theorem, together with the fact that all decomposition factors have support on intersections of hyperplanes(Lemma \ref{Prop32}), gives the number of decomposition factors of  $\rm{M_{\alpha}^{\beta}}$ and their support.

\begin{thm}
\label{Thm1}Let $A=\{H_{1},...,H_{m}\}$ be a  normal crossing arrangement on $\bC^n$.
 Let $\rm{H_S}=\rm{H_{i_{i}}\cap...\cap H_{i_{k}}}$ be a flat. Consider decomposition factors of  $\rm{M_{\alpha}^{\beta}}$.
 \begin{itemize}
 \item[(i)] There is at most one decomposition factor with support on $H_S$.
\item[(ii)]There is exactly one decomposition factor with support on $H_S$ if and only if $\beta_{i_{1}},...,\beta_{i_{k}}\in\bZ.$
 \end{itemize}
 \end{thm}
\begin{proof}By Lemma \ref{reductionlemma} we may assume that  $\cap_{i=1}^m H_i$ is either a point or $\emptyset$. If this intersection is a point $p$ then as stated in Lemma \ref{lemma:unique}, the arrangement is central and $m\leq n$, and the theorem is easy to prove, see  \cite[Prop. 3.1]{TARB}. (The idea is that by a basis change it may be assumed that $\alpha_i=x_i,\  i=1,...,m$ are coordinates at $p$ , and  then $\rm{M_{\alpha}^{\beta}}$ is isomorphic to the external tensor product of $\rm{M_{\alpha_i}^{\beta_i}}\  i=1,...,m$ on $\bC^1$.) 

This gives the induction basis for an inductive proof of (i) and (ii) on the number of hyperplanes in the arrangement, and it also allows us to assume that the intersection of all the hyperplanes in our configurations is empty. Assume that (i)  and (ii) are true for all normal crossing arrangements with $m$ hyperplanes, and let $A=\rm{\{H_{i}\}_{i=1}^{m+1}}$ be a normal crossing arrangement with $m+1$ hyperplanes.
\begin{itemize}
\item[(i)]
For $H_S=\bC^n$ the statement of (i) follows trivially, since $\rm{M_{\alpha}^{\beta}}$ has rank 1 as a module over $O_U$, where $U=\bC^n\setminus V(\alpha)$. Let $H_S\neq \bC^n$ be a flat of $A$. By Lemma \ref{lemma:unique} there is a hyperplane, which we may assume is $H_{m+1}$, such that $\rm{H_S\nsubseteq H_{m+1}}$. Let $\rm{U=X\setminus H_{m+1}}.$ By Proposition~\ref{Prop21} any decomposition factor $\rm{F}$ of $\rm{M_{\alpha}^{\beta}}$ with $\rm{SuppF=H}$, satisfies that the restriction $\rm{F|_{U}}$ is a non-zero irreducible $D_U$-module. By the proof of Lemma \ref{Prop33} the $D_U$-module $\rm{M_{\alpha}^{\beta}}\vert_U$ is isomorphic to the module $\rm{M_{\tilde \alpha}^{\tilde \beta}}\vert_U$ twisted by an automorphism, where $\tilde{\alpha}=\alpha_{1}...\alpha_{m}$ and $\tilde{\beta}=(\beta_{1},...,\beta_{m})$. The twisting does not change the number of decomposition factors nor their support. By the induction assumption and Proposition~\ref{Prop21} $\rm{M_{\tilde \alpha}^{\tilde \beta}}\vert_U$ has at most one decomposition factor on $H_S\cap U$. So the same is true of $\rm{M_{\alpha}^{\beta}}$ and by induction (i) is proved.

\item[(ii)] Similarily, we can use induction to prove one direction of (ii). Assume that (ii) is true for all normal crossing configurations with m hyperplanes, and that we have an arrangement $\rm{\{H_{i}\}_{i=1}^{m+1}}$. Let $\rm{H= H_{j_{1}}\cap...\cap H_{j_{s}}}$ be a flat such that  $\beta_{j_1}\in\bC\setminus\bZ$. By Lemma~\ref{blem} there exists $i$ such that $\rm{H\nsubseteq H_{i}}$.  By $(i)$ there exists at most one decomposition factor of $\rm{M_{\alpha}^{\beta}}$ with support H. Assume that there exists a decomposition factor of $\rm{M_{\alpha}^{\beta}}$, say F, with $\rm{SuppF=H}.$  Let $\rm{U=\bC^{n}\setminus H_{i}}.$ Then $\rm{F|_{U}}$ is simple  and $\rm{SuppF|_{U}=H\cap U.}$ Let $\tilde{\alpha}=\alpha_{1}...\hat \alpha_{i}...\alpha_{m+1}$ and $\tilde{\beta}=(\beta_{1},...,\hat\beta_{i},...,\beta_{m+1})$. Again arguing as in the proof of (i), we get that $\rm{M_{\tilde \alpha}^{\tilde \beta}}$ has a decomposition factor with support on $H$. But this is a contradiction to the induction assumption. Therefore there is no decomposition factor of $\rm{M_{\alpha}^{\beta}}$ which has support on H.

The converse of (ii) may be proven directly. Without loss of generality assume that $\rm{H=\cap_{i=1}^{k}H_{i}}$ such that $\beta_{1},...,\beta_{k}\in\bZ.$ Since we have a normal crossings arrangement, we may, by affine base change, assume that $\alpha_i=x_i,\ i=1,...,k$. Then $\rm{M_{\alpha}^{\beta}}\cong\bC[x]_{\alpha}\tilde{\alpha}^{\tilde{\beta}},$ where $\tilde{\alpha}^{\tilde{\beta}}=\alpha_{k+1}^{\beta_{k+1}}...\alpha_{m}^{\beta_{m}}.$ We first prove that there is one subfactor of $\rm{M_{\alpha}^{\beta}}$ with support on $\rm{H}$.
Let $\rm{N}\subset \bC[x]_{\alpha}$ be the vector space over $\bC$ generated by the set  
$$
\{\prod_{j=1}^n x_j^{r_j}\prod_{l=k+1}^{m}\alpha_{l}^{s_{l}}\tilde{\alpha}^{\tilde{\beta}}:{r_{j}\geq0\ \mathrm{if}\ j=k+1,...,n}\}.
$$  As a $\bC[x]$-module this is isomorphic to the localization of $\bC[x]_{x_1...x_k}$ to the open subset $\bC^n\setminus \cup_{j=k+1}^mH_j$. It is clearly an $A_n$-module.
Further 
let $\rm{N^{+}}\subset \rm{N}$ be the vector space over $\bC$ generated by the set 
$$\{\prod_{i=1}^n x_j^{r_j}\prod_{l=k+1}^{m}\alpha_{l}^{s_{l}}\tilde{\alpha}^{\tilde{\beta}}\colon r_{k+1}\geq0,...,r_{n}\geq0\ \&\
\exists l\colon \ 1\leq l\leq k \  \&\ r_l\geq 0\}$$
As a $\bC[x]$-module this is isomorphic to the localization of $\Sigma_{j=1}^k\bC[x]_{x_1...\hat x_j...x_k}$ to $\bC^n\setminus \cup_{j=k+1}^mH_j$. Again it is clearly an $A_n$-module, and $\rm{N/N^{+}}$ is  nonzero, since the quotient 
$$
\bC[x]_{x_1...x_k}/(\Sigma_{j=1}^k\bC[x]_{x_1...\hat x_j...x_k})
$$
is non-zero, with support on $H$, and $H\cap\bC^n\setminus \cup_{j=k+1}^mH_j\neq \emptyset$, by Lemma \ref{lemma:unique} (ii).
Then  $\rm{N/N^{+}}$ is a nonzero $A_{n}-$module, which is a decomposition factor of $\rm{M_{\alpha}^{\beta}}$.  As noted above $\rm{Supp(N/N^{+})=V(\alpha_{1},...,\alpha_{k})=\cap_{i=1}^{k}H_{i}=H}, $  and hence there must be at least one irreducible decomposition factor with support on $H$. \end{itemize}
This completes the proof of the theorem.
\end{proof}
It should be noted that the proof only uses that $\alpha_i$ such that $\alpha_i(p)=0$ form part of a system of parameters at each point $p$, and hence works in this situation for an arbitrary smooth affine variety.

\begin{definition}
Let A be a hyperplane configuration in $\bC^{n}$ and $\mathcal{H}$ be the set of all
nonempty intersections of hyperplanes in A, including $\bC^{n}$ itself considered as the intersection
over the empty set. Define a relation $x \leq y$ in $\mathcal{H}$ if $x\supseteq y$ (as subsets of $\bC^{n}$ ). In other words, $\mathcal{H}$ is partially ordered by reverse inclusion. We call $\mathcal{H}$ the intersection poset
of A.
\end{definition}
We can not give formulas for the cardinality $|\mathcal{H}|$ valid for all normal crossings arrangements, but in the special case when A is a hyperplane configuration in general position this is possible. In that case any intersection of less than or equal to $n$ hyperplanes gives a flat ( see \cite{TeraoOda}) :
$$|\mathcal{H}|=\sum_{k=0}^{n}{m\choose{k}}.$$
 
\begin{corollary}\label{Thm3.4}Let $\rm{M_{\alpha}^{\beta}=\bC[x]_{\alpha}{\alpha}^{\beta}},$ where $\alpha^{\beta}=\alpha_{1}^{\beta_{1}}...\alpha_{m}^{\beta_{m}}$ and the set $\{\alpha_{1},...,\alpha_{m}\}$ defines a hyperplane configuration in general position. Assume that $\beta_{1},...,\beta_{k}\in\bZ$ and $\beta_{k+1},...,\beta_{m}\in\bC\setminus\bZ.$ Then the number of decomposition facors of $\rm{M_{\alpha}^{\beta}}$, $c(\rm{M_{\alpha}^{\beta}})$ is given by: $$\rm{c(M}_{\alpha}^{\beta})=\sum_{j=0}^{n}{k\choose j},$$ where ${k\choose j}=0$ if $j>k.$ 
\end{corollary}
\begin{proof} In Theorem~\ref{Thm1} we proved that $\rm{M_{\alpha}^{\beta}}$  has exactly one  decomposition factor for each flat $\rm{H}$, where $\rm{H=H_{i_{1}}\cap...\cap H_{i_{s}}}$ and $\beta_{i_{1}},...,\beta_{i_{s}}\in\bZ$. The sum in the theorem clearly counts the number of those, keeping in mind that by the hypothesis on the arrangement, each non-empty subset $S\subset\{1,...,m\}$ corresponds to a unique flat.
\end{proof}

\section{Example}
\subsection{The $\rm{A_{2}}$-module $\rm{M_{\alpha}^{\beta}=\bC[x,y]_{xy(x+y+1)}x^{\beta_{1}}y^{\beta_{2}}(x+y+1)^{\beta_{3}}}$.} We will describe the decomposition factors of $\rm{M_{\alpha}^{\beta}}$ for different cases of $\beta_{1},\beta_{2}$ and $\beta_{3}$, as an example for the previous section. The following Theorem summarizes the results.
\unitlength=1mm
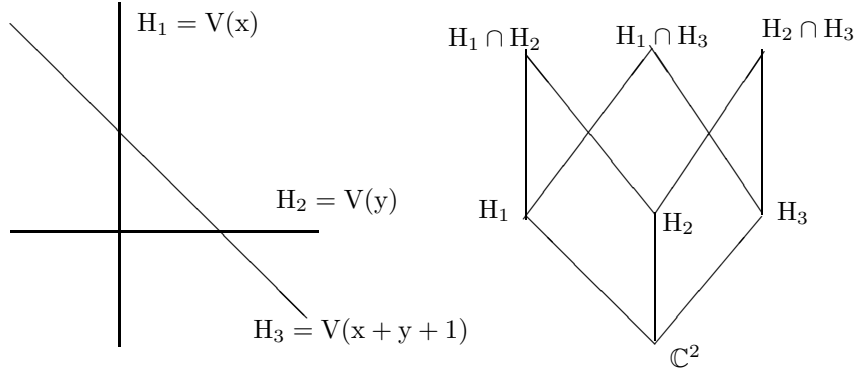
\begin{figure}
\centering
\begin{picture}(121.92,58.88)(12.48,42.56)
\put(38.08,95.84){\line(0,-1){45.6}}
\put(23.68,65.52){\line(1,0){40.52}}
\put(23.68,93.12){\line(1,-1){39.04}}
\put(108.48,51.10){\line(0,1){16.90}}
\put(108.48,50.66){\line(-1,1){16.96}}
\put(108.48,50.65){\line(5,6){14.08}}
\put(91.52,67.2){\line(0,1){22.4}}
\put(108.48,67.84){\line(-4,5){16.96}}
\put(122.56,67.84){\line(0,1){21.76}}
\put(108.48,67.84){\line(2,3){14.4}}
\put(91.2,67.2){\line(3,4){16.96}}
\put(122.56,68.16){\line(-2,3){14.72}}
\put(55.68,50.56){\makebox(16.32,6.4)[lb]{$\rm{H_{3}=V(x+y+1)}$}}
\put(40.4,91.2){\makebox(18.56,8.96)[lb]{$\rm{H_{1}=V(x)}$}}
\put(58.68,67.84){\makebox(17.92,6.72)[lb]{$\rm{H_{2}=V(y)}$}}
\put(110.48,47.56){\makebox(10.56,8)[lb]{$\bC^{2}$}}
\put(85.28,66.88){\makebox(10.24,7.04)[lb]{$\rm{H_{1}}$}}
\put(109.5,65.12){\makebox(7.04,6.72)[lb]{$\rm{H_{2}}$}}
\put(124.56,66.48){\makebox(8.96,8.32)[lb]{$\rm{H_{3}}$}}
\put(81.28,89.6){\makebox(10.24,10.56)[lb]{$\rm{H_{1}\cap H_{2}}$}}
\put(103.36,89.92){\makebox(11.84,10.56)[lb]{$\rm{H_{1}\cap H_{3}}$}}
\put(122.56,90.88){\makebox(11.84,10.56)[lb]{$\rm{H_{2}\cap H_{3}}$}}
\end{picture}
\caption{A hyperplane configuration defined by the polynomial $\alpha=xy(x+y+1)$ and the Hasse diagram of its intersection poset.}
\label{fig:pic}
\end{figure}
\begin{thm} Consider the $A_{2}-$module  $\rm{M_{\alpha}^{\beta}=\bC[x,y]_{xy(x+y+1)}x^{\beta_{1}}y^{\beta_{2}}(x+y+1)^{\beta_{3}}}$. 
\begin{itemize}
\item[(i)] If $\beta_{1},\beta_{2},\beta_{3}\in \bC\setminus\bZ$, then $c(\rm{M_{\alpha}^{\beta}})=1$ and hence $\rm{M_{\alpha}^{\beta}}$ is irreducible.
\item[(ii)] If $\beta_{1}\in \bZ$ and $\beta_{2},\beta_{3}\in \bC\setminus\bZ,$ then $c(\rm{M_{\alpha}^{\beta}})=2.$
\item[(iii)] If $\beta_{1},\beta_{2}\in \bZ$ and $\beta_{3}\in \bC\setminus\bZ,$ then $c(\rm{M_{\alpha}^{\beta}})=4.$
\item[(iv)] If $\beta_{1},\beta_{2},\beta_{3}\in \bZ$, then $c(\rm{M_{\alpha}^{\beta}})=7.$
\end{itemize}
\end{thm}
\begin{proof}The linear functions $x,y$ and $x+y+1$ define a general position hyperplane configuration in $\bC^{2}.$ Therefore, by the previous section, the number of the decomposition factors of $\rm{M_{\alpha}^{\beta}}$ is the number of flats $\rm{H}$ such that the corresponding $\beta's$ of the linear functions of the hyperplanes defining $\rm{H}$ are integers.
\begin{itemize}
\item[(i)]If $\beta_{1},\beta_{2},\beta_{3}\in \bC\setminus\bZ,$ then $\rm{H=\bC^{2}},$ considered as the empty intersection, is the only flat such that the $\beta's$ of the linear functions defining the hyperplanes defining $\bC^{2}$ are integers. Therefore $\rm{M_{\alpha}^{\beta}}$ has only one decomposition factor with support on $\bC^{2}$ and hence $\rm{M_{\alpha}^{\beta}}$ is irreducible with support on $\bC^{2}.$
\item[(ii)] If $\beta_{1}\in \bZ$ and $\beta_{2},\beta_{3}\in \bC\setminus\bZ,$ then $\rm{H_{1}=\bC^{2}},$ considered as the empty intersection and $\rm{H_{2}=V(x)}$ are the flats with corresponding integer $\beta's$. We can also see this by the following composition series:
$$\bC[x,y]_{y(x+y+1)}y^{\beta_{2}}(x+y+1)^{\beta_{3}}\subset\bC[x,y]_{xy(x+y+1)}x^{\beta_{1}}y^{\beta_{2}}(x+y+1)^{\beta_{3}}$$ 
where
$$\bC[x,y]_{y(x+y+1)}y^{\beta_{2}}(x+y+1)^{\beta_{3}}$$ and $$\bC[x,y]_{xy(x+y+1)}x^{\beta_{1}}y^{\beta_{2}}(x+y+1)^{\beta_{3}}/\bC[x,y]_{y(x+y+1)}y^{\beta_{2}}(x+y+1)^{\beta_{3}}$$ are easily seen to be irreducible $A_{2}-$modules. 
Therefore, as predicted, $\rm{c(M_{\alpha}^{\beta})=2}$, with one decomposition factor with support on $\rm{V(x)}$ and one decomposition factor with support on the whole space $\bC^{2}.$ 
\item[(iii)]If $\beta_{1},\beta_{2}\in \bZ$ and $\beta_{3}\in \bC\setminus\bZ,$ then $\rm{H_{1}=\bC^{2}},$ considered as the empty intersection and $\rm{H_{2}=V(x),H_{3}=V(y)}$ and $\rm{H_{4}=V(x,y)}$ are the flats with corresponding integer $\beta's$.  Consider the following sequence of vector spaces over $\bC$: $$R_{0}\subset R_{1}\subset \rm{M_{\alpha}^{\beta}},$$ where  $R_{0}=\bC[x,y]_{(x+y+1)}(x+y+1)^{\beta_{3}},$ and $R_{1}$ is generated by the set $$\{\frac{\alpha_{j}^{r_{j}}(x+y+1)^{s}}{\alpha_{i}^{s_{i}}(x+y+1)^{\beta_{3}}}:s_{j}\geq0,s_{i}\geq1,\ \rm{for} \ i,j=1,2\}.$$Therefore $\rm{c(M_{\alpha}^{\beta})=4}$ with one decomposition factor with support on each $\rm{H_{2},H_{3},H_{4}}$, one decomposition factor with support on $\bC^{2}.$ 
\item[(iv)]If $\beta_{1},\beta_{2},\beta_{3}\in \bZ,$ then $\rm{H_{1}=\bC^{2}},$ considered as the empty intersection and $\rm{H_{2}=V(x),H_{3}=V(y),H_{4}=V(x+y+1),H_{5}=V(x,y),H_{6}=V(x,x+y+1)}$ and $\rm{H_{7}=V(y,x+y+1)}$ are the flats with their $\beta's$ of the linear functions defining the hyperplanes defining $\rm{H_{1},...,H_{7}}$ are integers. Therefore $\rm{c(M_{\alpha}^{\beta})=7}$, with exactly one decomposition factor with support on each of the flats. 
\end{itemize} 
\end{proof}

\section{Resolution of singularities} Let $X=\bC^n\setminus V(\alpha)$, and $\pi\colon \tilde \bC^n\to \bC^n$ be a resolution of the singularities of $V(\alpha)$. Let $\rm{E=\bigcup_{i=1}^{s} E_{i}}$, where $E_{i}$ are irreducible, be the exceptional divisor, and $Z=\pi(E)$ the center, so that $\tilde \bC^n\setminus E\cong \bC^n\setminus Z$. We then have that 
\begin{equation}
\label{eq:div}
{\rm Div}\pi^{*}(\alpha_{i})=\tilde{H_{i}}+r_{1}^{i}E_{1}+...+r_{s}^{i}E_{s},
\end{equation}
 where $\tilde{H_{i}}$ is the proper transform of $H_{i}=V(\alpha_{i})$ and $r_j^i\in \bZ_+$. 
 Hence in an open affine $U\subset  \tilde \bC^n$,
 \begin{equation}
 \label{eq:res}
 \pi^{*}(\alpha_{i})|_{U}=(\tilde{\alpha}_{i}^{\rm{U}}(\gamma_{1}^{\rm{U}})^{r_{1}^{i}}...(\gamma_{s}^{\rm{U}})^{r_{s}^{i}})\subset \mathcal{O}_{U}
\end{equation}
  where $\tilde{H_{i}}\cap U=V(\tilde{\alpha}^U_{i})$ and $\rm{U\cap E_{j}=V(\gamma_{j}^{U})}.$ Some of the intersections $\tilde{H_{i}}\cap U$ or $\rm{U\cap E_{j}}$ may be empty, but we supress this from notation, so as to not make it cumbersome. The properties  of $\tilde \bC^n$ that we need are
\begin{itemize}
 \item[i)] that there is an 
affine cover $\{\rm{ U_{i}\subset\tilde\bC^{n}},i=1,...,r\}$, where each $U=U_i$ is isomorphic to an open affine subvariety of finite type of $ \bC^n$ 
\item[ii)] $\pi^{-1}V(\alpha)$ is a normal crossing divisor. 

\end{itemize}
For these standard properties see e.g. \cite{HJ1}, or the construction in the present case in \cite{CDCP2}.
We will use the preceding section to study the D-module pullback  $\pi^*M_\alpha^\beta $. The last statement above implies that Theorem \ref{Thm1} applies to the restriction of this module to each $U_i$. 

\begin{lemma} Let $U=U_i,\ i=1,...,r.$ Then
$\pi^*M_\alpha^\beta \vert_{U}\cong M_{\tilde \alpha}^{\tilde \beta} $  as $D_{U}$-modules, where  
$\tilde \alpha=\tilde{\alpha}_{1}^{\rm{U}}....\tilde{\alpha}_{m}^{\rm{U}}\gamma_{1}^{\rm{U}}..\gamma_{s}^{\rm{U}}$.(In the product we only use the factors for which the corresponding intersections $\tilde{H_{i}}\cap U$ or $\rm{U\cap E_{j}}$ are non-empty.) 
Further (with the same restriction)
 $$(\tilde\beta_1,\tilde\beta_2,...,\tilde\beta_{m+s})=(\beta_1,\beta_2,...,\beta_m, \sum_{i=1}^{m}r_{1}^{i}\beta_{i},...,\sum_{i=1}^{m}r_{s}^{i}\beta_{i}).$$
\end{lemma}
Using the last section, we immediately obtain the decomposition factors of $\pi^*M_\alpha^\beta$.
Extend the notation by defining $ \tilde{\rm H}_{m+i}:={\rm E}_i,\ i=1,...,s.$

\begin{corollary}\label{Cor55} The module $\pi^{*}(\rm{M_{\alpha}^{\beta}})$ has exactly one decomposition factor for each flat $\rm{\tilde{H}}={\rm {\tilde H}}_{i_{1}}\cap...\cap{\rm{ \tilde {H}}}_{i_{k}}$ such that $\tilde{\beta}_{i_{1}},...,\tilde{\beta}_{i_{k}}\in\bZ$. 
\end{corollary}
\begin{proof}Since the $U_i$ cover $\tilde \bC^n$, there is to any decomposition factor $N$ of $\pi^{*}(\rm{M_{\alpha}^{\beta}})$ a $U_i$ such that $N\vert_{U_i}\neq 0$. This restriction has, by Proposition \ref{Prop32}, support on an intersection 
$\rm{\tilde{H}}\cap U_i=\tilde{H}_{i_{1}}\cap...\cap \tilde{H}_{i_{k}}\cap U_i$, 
and by  Theorem   \ref{Thm1} and the preceding Lemma, this implies that $\tilde{\beta}_{i_{1}},...,\tilde{\beta}_{i_{k}}\in\bZ$. Furthermore there is by  Theorem  \ref{Thm1} at most one. Conversely,
assume that $\rm{\tilde{H}}=\tilde{H}_{i_{1}}\cap...\cap \tilde{H}_{i_{k}}$ is such that $\tilde{\beta}_{i_{1}},...,\tilde{\beta}_{i_{k}}\in\bZ$. There is $U_i$ such that $\rm{\tilde{H}}\cap U_i\neq \emptyset $, and hence there is by  Theorem  \ref{Thm1} a decomposition factor of $\pi^{*}(\rm{M_{\alpha}^{\beta}})\vert _{U_i}$  with support on $\rm{\tilde{H}}\cap U_i$. By Corollary \ref{cor:open}
this means that there is a decomposition factor of $\pi^{*}(\rm{M_{\alpha}^{\beta}})$ with support on  $\rm{\tilde{H}}$.
\end{proof}
\begin{corollary}\label{Cor56} The module  $\pi^{*}(\rm{M_{\alpha}^{\beta}})$ is irreducible if and only if $$\beta_{1},...,\beta_{m},\sum_{i=1}^{m}r_{i}^{1}\beta_{i},...,\sum_{i=1}^{m}r_{i}^{s}\beta_{i}\in \bC\setminus\bZ.$$
\end{corollary}
We can use this and the decomposition theorem to get a sufficient criterion for when $\rm{M_{\alpha}^{\beta}}$ is irreducible. We will use the decomposition theorem in the general form stated by Kashiwara \cite{Kashiwara2} and recently proved by Mochizuki \cite{Mochizuki}. (Our module is regular, so our application actually could have referred to old versions of the decomposition theorem.) 
\begin{thm}
\label{Thm:irr}
Assume that $r^i_j$ are defined by the property (\ref{eq:div}) of a resolution of the singularities of $V(\alpha)$. Then $\rm{M_{\alpha}^{\beta}}$ is irreducible if  $\beta_{1},...,\beta_{m},\sum_{i=1}^{m}r_{i}^{1}\beta_{i},$ $...,\sum_{i=1}^{m}r_{i}^{s}\beta_{i}\in \bC\setminus\bZ$. 
\end{thm}
\begin{proof}Given the condition in the theorem, by Corollary~\ref{Cor56}, $\rm{\pi^{*}(M_{\alpha}^{\beta})}$ is irreducible. Consider now the following pullback diagram:
\begin{equation}\label{Dgm1}
\xymatrix{
\tilde{U}\ar[r]^{ j'}\ar[d]_{\pi_{|_{\tilde U}}}&  
\tilde{\bC^{n}}\ar[d]_{\pi} \\
U\ar[r]^{j}& \bC^{n}} 
\end{equation}
where $U=\bC^{n}\setminus V(\alpha)$ and $j'$ and $j$ are the respective inclusions. By definition  $\rm{M_{\alpha}^{\beta}=j_{*}({M_{\alpha}^{\beta}}_{|_{U}})}$. Since $j$ is an affine map, this remains true interpreted in the derived category $D(D_{\bC^n})$, where $j_*$ now is the derived functor of $\Oc_{\bC^n}-$modules(which is also the direct D-module image $j_+$). Similarily $j'_*$ and
${\pi^*_{|_{\tilde U}}}$ are exact functors, the latter since $\pi_{|_{\tilde U}}$  is an isomorphism.
By Theorem 1.7.3 of \cite{Hotta}, pullbacks behave nicely in the derived category $D(D_{\tilde \bC^n})$, so
\begin{equation}
\label{eq:pullbackdescr}
\rm{\pi^*M_{\alpha}^{\beta}}=\rm{\pi^*j_{*}({M_{\alpha}^{\beta}}_{|_{U}})=j'_{*}\pi_{|_{\tilde U}}^{*}({M_{\alpha}^{\beta}}_{|_{U}})},
\end{equation}
and consequently 
\begin{equation}\label{eq:module}
H^r\pi^*\rm{M_{\alpha}^{\beta}}=0,\ if \ r\neq 0.
\end{equation} 
Now, in the derived category,
\begin{equation}
\label{eq:nohighercohomology}
\rm{\pi_{+}\pi^{*}(M_{\alpha}^{\beta})=\pi_{+}j'_{*}(\pi_{|_{\tilde U}}^{*}({M_{\alpha}^{\beta}}_{|_{U}}))}=j_{*}(\pi_{|_{\tilde U}})_{+}\pi_{|_{\tilde U}}^{*}({M_{\alpha}^{\beta}}_{|_{U}}))=j_{*}({M_{\alpha}^{\beta}}_{|_{U}}))=M_{\alpha}^{\beta}.
\end{equation} 
But by the Decomposition Theorem (see e.g.  \cite{Mochizuki}), $H^0\rm{\pi_{+}\pi^{*}(M_{\alpha}^{\beta})}$ is semisimple. However $\rm{M_{\alpha}^{\beta}}$ is indecomposible, since it is rank 1 over the generic point and torsion free as a $\bC[x]-$module, and hence it is irreducible.


\end{proof}

\section{Plane case} 
In this section we will exemplify the preceding, by describing the pullback to a resolution of singularities of the $A_{2}-$module $\rm{M_{\alpha}^{\beta}=\bC[x,y]_{\alpha}\alpha^{\beta}}$, where $\alpha=\prod\limits_{i=1}^{m}\alpha_{i},$ $\alpha_{1}=x,$ $\alpha_{2}=y$ and $\alpha_{i}=x+c_{i}y$ for all $i=3,..,m$ and $ c_{i}\neq c_{j}$ for all $i\neq j$ and $\beta=(\beta_{1},\beta_{2},...,\beta_{m})\in \bC^m$. 
We only need to blow up the origin. This is the locus:
$$\tilde{\mathbb{C}^{2}}=\{(x,y),[W_{0},W_{1}]):xW_{1}=yW_{0}\}\subset \mathbb{C}^{2}\times \mathbb{P}^{1}$$ together with the map $$\pi:\tilde{\mathbb{C}^{2}}\longrightarrow\mathbb{C}^{2}$$ which is the restriction of the projection of $\mathbb{C}^{2}\times \mathbb{P}^{1}$ onto the first factor. Let $U_{1}\subset \tilde{\mathbb{C}^{2}}$ be the open subset given by $W_{0}\neq 0$. In terms of Euclidian coordinates, $w_{1}=\frac{W_{1}}{W_{0}}$ and we can write: $$U_{1}=\{(x,y),(w_{1})):xw_{1}=y\}=\{(x,xw_{1},w_{1})\}\subset \mathbb{C}^{2}\times \mathbb{C}^{1}.$$ From this description we see that $U_{1}\cong \mathbb{C}^{2}$ with coordinates $x,w_{1}$. The restriction $\pi|_{U_{1}},$ is given by: $$\pi(x,w_{1})=(x,xw_{1}).$$
Similarily let $U_{2}\subset \tilde{\mathbb{C}^{2}}$ be given by $W_{1}\neq 0$, and put  $w_{0}=\frac{W_{0}}{W_{1}}$. Then: $$U_{2}=\{(x,y),(w_{0})):x=yw_{0}\}=\{(yw_{0},y,w_{0})\}\subset \mathbb{C}^{2}\times \mathbb{C}^{1}.$$ So $U_{2}\cong \mathbb{C}^{2}$ with coordinates $y,w_{0}$, and the restriction $\pi|_{U_{2}},$ is given by: $$\pi(y,w_{0})=(yw_{0},y).$$ Clearly  $\tilde{\mathbb{C}^{2}}=U_{1}\cup U_{2}$ and hence $\{U_{1},U_{2}\}$ is an affine cover of $\tilde{\mathbb{C}^{2}}.$ 
The exceptional divisor $E$ is given in the chart $U_1$ by $x=0$, and in $U_2$ by $y=0$. Furthermore
in $U_1$
$$
\pi^*(\alpha_1)=(x), \pi^*(\alpha_2)=(xw_1),\pi^*(\alpha_i)=(x(1+c_iw_1)), \ i=3,...,m,
$$
and in
in $U_2$
$$
\pi^*(\alpha_1)=(yw_0), \pi^*(\alpha_2)=(y),\pi^*(\alpha_i)=(y(w_0+c_i)), \ i=3,...,m.
$$
Hence, $Div\pi^*(H_i)=\tilde H_i+E, \ i=1,...,m$, (and $U_1={\tilde \bC}^n\setminus H_1$ and $U_2={\tilde \bC}^n\setminus H_2$.) This gives us the numerical description of the resolution that we need to apply Corollary \ref{Cor55}.

 \begin{proposition}\label{bethm1} Consider $\rm{\pi_{2}^{*}(M_{\alpha}^{\beta})}$ . Let $I=\{i\vert \beta_i\in \bZ\}$ contain $k$ elements, and let $\vert \beta\vert=\sum_{i=1}^m\beta_i$ .
 \begin{itemize}
\item[(i)]If $\vert \beta\vert \in \bZ$,  then there is exactly one decomposition factor with support on each of $\bC^2$, $E$, $H_i, \ i\in I$ and one on each intersection $H_i\cap E,\ i\in I$, and these are all decomposition factors. All together $2(k+1)$.
\item[(ii)]If $\vert \beta\vert \in\bC\setminus  \bZ$,  then there is exactly one decomposition factor with support on $\bC^2$,  and one on each $H_i, \ i\in I$ only. All together $k+1$.
\end{itemize}
\end{proposition}

\unitlength=1mm
\begin{figure}
\centering
\begin{picture}(141.44,61.44)(29.44,39.36)
\put(61.44,100.48){\line(0,-1){60.8}}
\put(29.76,70.08){\line(1,0){61.44}}
\put(36.8,91.2){\line(6,-5){54.72}}
\put(84.16,93.44){\line(-1,-1){48.64}}
\put(50.24,100.16){\line(2,-5){24}}
\put(75.84,99.2){\line(-1,-2){29.12}}
\put(31.68,81.6){\line(5,-2){59.84}}
\put(88.96,83.84){\line(-2,-1){59.52}}
\put(139.84,100.8){\line(0,-1){55.44}}
\put(108.8,70.4){\line(1,0){62.08}}
\put(108.8,65.56){\line(1,0){62.08}}
\put(108.8,60.32){\line(1,0){62.08}}
\put(107.84,75.64){\line(1,0){62.08}}
\put(107.84,80.6){\line(1,0){62.4}}
\put(107.84,85.64){\line(1,0){62.4}}
\put(108.48,55.92){\line(1,0){62.4}}
\end{picture}
\caption{The hyperplane configuration before and after the blowup in one of the affine charts.}
\label{fig:blowup1}
\end{figure}
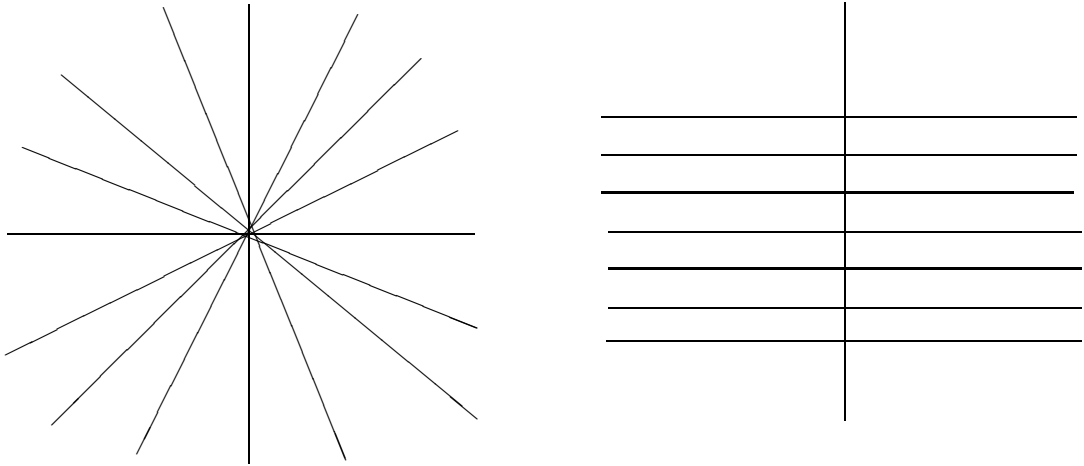

We can compare the criterion for irreducibility of $\rm{M_{\alpha}^{\beta}}$ that was obtained in Theorem \ref{Thm:irr}, with what was proved in the plane case in \cite{TARB} and see that in the plane it also gives a necessary criterion. We ignore whether this is true in general.
\begin{corollary} 
\label{finalthm}In the plane case, $\pi^{*}(\rm{M_{\alpha}^{\beta}})$ is irreducible if and only if $\rm{M_{\alpha}^{\beta}}$ is irreducible if and only if $\beta_{1},....,\beta_{m},\sum_{i=1}^{m}\beta_{i}\in\bC\setminus\bZ.$
\end{corollary}
\begin{proof}
See  \cite{TARB} .
\end{proof}
\section*{Acknowledgements}We gratefully acknowledge the encouragement we got while working on these topics from Demissu Gemeda, Addis Ababa University and would like to dedicate this note to his memory. We also want to thank Jan-Erik Bj\"ork and Rolf K\"allstr\"om for
interesting, useful and fun discussions. Finally the first author also acknowledges the economical support of ISP, Uppsala University.

\end{document}